\documentclass[12pt]{article}

\usepackage{amsmath}
\usepackage{amsthm}
\usepackage{amssymb}
\usepackage{amstext} 
\usepackage{float}
\usepackage{graphicx}
\usepackage{subfig}
\usepackage{bbm}
\usepackage[pagewise]{lineno}

\newtheorem{theorem}{Theorem}

\newtheorem{rem}[theorem]{Remark}

\usepackage{listings}
\usepackage{nicefrac}
\usepackage{color}
\usepackage{epstopdf}

\usepackage{amssymb}

\usepackage{lineno}

\usepackage{tabularx}
\newcolumntype{+}{>{\global \let \currentrowstyle \relax}}
\newcolumntype{^}{>{\currentrowstyle }}

\usepackage{booktabs}
\usepackage{multirow}

\usepackage{lscape}

\allowdisplaybreaks

\newcommand\BibTeX{{\rmfamily B\kern-.05em \textsc{i\kern-.025em b}\kern-.08em
		T\kern-.1667em\lower.7ex\hbox{E}\kern-.125emX}}

\begin{document}

\title{Optimal Control for Tuberculosis with Exogenous Reinfection and Stigmatization}

\author{Remilou F. Liguarda\\
	{\small MSU-IIT Iligan, Andres Bonifacio Avenue, Tibanga},\\ 
	{\small	9200 Iligan City, Philippines}\\
	{\small E-Mail:remilou.liguarda@g.msuiit.edu.ph}\\[.3cm]
	Wolfgang Bock\\
	{\small Technomathematics Group}\\
	{\small	University of Kaiserslautern}\\
	{\small	P.\ O.\ Box 3049, 67653 Kaiserslautern, Germany}\\
	{\small E-Mail:bock@mathemaik.uni-kl.de}\\[.3cm]
	Randy L. Caga-anan\\
	{\small MSU-IIT Iligan, Andres Bonifacio Avenue, Tibanga},\\ 
	{\small	9200 Iligan City, Philippines}\\
	{\small E-Mail:randy.caga-anan@g.msuiit.edu.ph}}

\maketitle

\begin{abstract}
	The effect of stigmatization is hindering the control of diseases. Especially in the case of exogenous reinfection, this effect can play a massive role. We develop in this paper, based on a tuberculosis model of Feng et.al. a model with exogenous reinfection and stigmatization. As in the base model, in presence of exogenous reinfection there exists multiple endemic equilibria. We solve an optimal control problem for a case scenario and show that both types of measures, those fighting stigmatization but also direct disease control have to be invoked.
\end{abstract}

\section{Introduction}

According to World Health Organization's Global Tuberculosis Report 2019 \cite{4}, TB was one of the top 10 causes of death worldwide in 2018 and it is also the leading killer of people with HIV and a major cause of deaths related to antimicrobial resistance. In 2018, a total of $1.5$ million people died from TB were recorded among which about $251 \ 000$ people who were also infected with HIV. In the same year, it was estimated globally that there were 10 million new cases of TB out of which $9\%$ people were co-infected with HIV. Moreover, the countries such as India, China, Indonesia, Philippines, Pakistan, Nigeria, Bangladesh, and South Africa were accounted for $66\%$ of the new cases.

TB is a treatable and curable disease. With an early diagnosis and treatment with antibiotics, most people who develop TB can be cured and onward transmission curtailed \cite{4}. Treatment will
only be effective if the patient completes the therapy which includes a combination of drugs recommended by the physicians \cite{7}. Poor compliance contributes to the worsening of the TB situation by increasing incidence and initiating drug resistance \cite{7}. The emergence of drug resistant TB, especially multidrug-resistant TB (MDR-TB) and extensively drug-resistant TB (XDR-TB), poses a substantial threat to TB control programs worldwide \cite{9}. In 2018, MDR-TB remains a public health crisis and a health security threat and among of its cases, 6.2\% were estimated to have XDR-TB \cite{3}. Treatment success rate at 56\%, remains low globally \cite{4}. Some factors reported to have a significant effect on adherence are: poor socioeconomic status, poverty, illiteracy, low level of education, unemployment, lack of effective social support networks, unstable living conditions, long distance from treatment centre, high cost of transport, high cost of medication, changing environmental situations, culture and lay beliefs about illness and treatment, and family dysfunction \cite{10}.

It is known that social norms and poverty lead to a stigmatization of tuberculosis (TB), i.e.~ patients with severe symptoms are not seeking medical treatment, hence stigmatization is thought to hinder TB control. \cite{CW}
In particular the social status referring to institutional and social factors affect the individual access of medical services \cite{WHO_soc, HvdM}. Stigma is a process that begins particular groups are characterized as being undesirable or disvalued, leading to  shame, disgust, and guilt see \cite{Link,Goffman}.  A good overview about stigmatization with tuberculosis can be found in \cite{CW}.

The dynamics of communicable diseases can be explored affectively with the help of mathematical models and provide useful information to the spread and control of the communicable diseases \cite{1}. For the case of TB a latent phase is important, see e.g.~\cite{FCC01}. There are several models based on different aspects of the disease as multistrain models \cite{Jung02, Feng02}, treatment and optimal control \cite{CC97, CC98, Ain17}, the modeling of HIV co-infection \cite{Ro09,13} and migration \cite{Zhou} just to name a few.

In this paper, a compartmental model for tuberculosis disease transmission considering stigmatization and enhanced reinfection is formulated. Particularly to the model from Feng et.al.~\cite{13}  we add the stigmatization effect by deviding  the active infected compartment into two, namely, those willing to seek medical treatmentand those not due to stigma. 
The basic reproductive number $R_0$ and equilibrium points of the model are determined and the stability analysis is performed considering the $R_0$. We show the existence of multiple endemic equilibrium points as in \cite{13}. We study optimal control strategies where we also include a control parameter which brings people not seeking treatment back to seek medical treatment. Our numerical results show the most effective strategy uses several different controls.

\section{Mathematical Model}

We took the model from Feng et.al.~\cite{13} as our base model, which describes the transmission of TB with exogenous reinfection. To add the stigmatization effect, we further divide the active infected compartment into two, namely, those willing to seek medical treatment $(I_S)$ and those not willing to seek medical treatment $(I_N)$ because of the stigma associated with being a TB patient. The entire population is classified into five classes: susceptible $(S)$, exposed $(E)$, infectious willing to seek medical attention $(I_S)$ and infectious but not willing to be treated $(I_N)$ and those under treatment or already treated $(T)$. We denote by $N$ the total population $S+E+I_S+I_N+T$ and $I=I_S+I_N$.

The model is governed by the following system of ordinary differential equations:
{
	\begin{eqnarray}
		\frac{d S}{d t} &=&\Lambda-\beta cS\frac{I}{N}-\mu S \label{eq1}\\
		\frac{d E}{d t} &=& \beta cS\frac{I}{N}-p\beta cE\frac{I}{N}-(\mu +k)E+\sigma \beta c T\frac{I}{N} \\
		\frac{d I_S}{d t} &=& \alpha p\beta c E \frac{I}{N}+\alpha kE-(\mu +r+d)I_S\\
		\frac{d I_N}{d t} &=& (1-\alpha)p\beta cE\frac{I}{N} +(1-\alpha)kE-(\mu +d)I_N\\
		\frac{d T}{d t} &=& rI_S-\sigma \beta cT\frac{I}{N}-\mu T \label{eq5}
	\end{eqnarray}
}

\begin{rem}
	The parameter $p$ represents in both models the one from Feng et. al. \cite{13} and the model (1)-(5) the role of reinfection. For $p=0$ we obtain for $\alpha \neq 1$ a TB model as Castillo-Chavez and Feng \cite{CCF97} with stigmatization. For $p \leq 1$ a reinfection is less likely than an infection. As pointed out in Feng et. al.\cite{13} $p>1$ maybe reasonable for HIV-infected individuals. 
\end{rem}

The parameters are described in Table \ref{tab:parameters}. These are the same parameters in \cite{13}, with the addition of the parameter $\alpha$ with values in the interval $[0,1]$, controlling the level of stigmatization. An $\alpha=1$ means the infectious infected are all seeking medical treatment where $\alpha=0$ means all of them avoid medical treatment. 
\begin{table}[htbp]
	\centering
	\begin{tabular}{c|c|c}
		\hline 
		Parameter	& Description & Value\\ 
		\hline 
		$\Lambda$ & recruitment rate &    588 humans/year \\
		\hline 
		$\beta c$ & transmission rate  from $S$ to $E$  &   2 /year \\ \hline 
		$\sigma$ & reduction of reinfection rate from $T$  &   0.9 (dimensionless) \\
		\hline 
		$\mu$ &  natural death rate &  0.0235 /year \\ 
		\hline 
		$k$ & progression rate from $E$ to $I$  &  0.0294 /year  \\ 
		\hline
		$d$ &  disease-induced death rate  & 0.05 /year  \\ 
		\hline
		$r$ & treatment rate  &  0.2906 /year \\ 
		\hline 
		$p$ & level of exogenous reinfection  &  0.4 (dimensionless)\\ 
		\hline 
		$\alpha$ & level of stigmatization  &  [0,1] (dimensionless)\\ 
		\hline 
	\end{tabular} 
	\caption{Model Parameters}
	\label{tab:parameters}
\end{table}

\section{Analysis of the system and equilibria}
The model will be analyzed in the biologically-feasible region, $\Gamma \subset \mathbb{R}_{+}^{5}$ with,
$$\Gamma=\{(S(t), E(t), I_S(t), I_N(t), T(t)) \in \mathbb{R}_{+}^{5}: 0 \leq S,E,I_S,I_N,T\}.
$$
Define $f: \mathbb{R}_{+}^{5} \rightarrow \mathbb{R}^{5}$ by $f(y)=\left(f_{1}(y), f_{2}(y), f_{3}(y), f_{4}(y),f_{5}(y)\right),$ where $f_{i}: \mathbb{R}_{+}^{5} \rightarrow \mathbb{R}^{5}$ is defined by
$f_{i}(y)=\frac{d y_{i}}{d t}$ for each with $i=1,2,3,4,5.$ Then the system $(1)-(5)$ allows us to define an initial value problem (IVP)
$$
\begin{array}{rcl}
\dot{\mathbf{y}} &=& f(y) \\
y\left(t_{0}\right) &=& y_{t_{0}}
\end{array}
$$
where $y=(S, E, I_S, I_N, T) .$ As all functions on the right hand side (R.H.S) of the model above are continuously differentiable on $\mathbb{R}_{+}^{5}$ if $y_{0} \in \mathbb{R}_{+}^{5}$ then IVP possesses a locally unique solution by Picard-Lindel\"of.

\begin{theorem}
	Let $y_{0} \in \mathbb{R}_{+}^{5} .$ Any solution $y(t)$ of IVP through $y_{0}$ is defined for all $t \geq 0,$ and the region $\Gamma \subset \mathbb{R}_{+}^{5}$ is positively invariant.
\end{theorem}

\noindent{\textbf{Proof:}}
Let $\left(S^{*}(t), E^{*}(t), I_S^{*}(t), I_N^{*}(t),T^{*}(t) \right)$ with $S^{*}(0)>0$, $E^{*}(0)>0$, $I_S^{*}(0)>0$, $I_N^{*}(0)>0$, $T^{*}(0)>0$ be a fixed solution of IVP (6)-(7) through $\left(S^{*}, E^{*}, I_S^{*}, I_N^{*}, T^{*}\right)$ on $[0, \xi),$ with $\xi$ real. With the continuity of the solution on
$[0, \xi),$ there exists $\delta$ with $0<\delta<\xi$ such that $S^{*}\left(t^{\prime}\right)>0$, $E^{*}\left(t^{\prime}\right)>0$, $I_S^{*}\left(t^{\prime}\right)>0$, $I_N^{*}\left(t^{\prime}\right)>0$, $T^{*}\left(t^{\prime}\right)>0$, for any $t^{\prime} \in[0, \delta]$,

$$\begin{aligned}
\frac{d S^{*}}{d t} &=\Lambda-\beta cS^{*}\frac{I^{*}}{N}-\mu S^{*} \\
& =\Lambda+S^{*}(-\beta c\frac{I^{*}}{N}-\mu)  \\
& \geq S^{*}(-\beta c\frac{I^{*}}{N}-\mu) \\
& \geq S^{*}(-\frac{\beta c}{N}-\mu) \\
\end{aligned}$$
since $N^{*}\geq I^{*}\geq 0 .$ 
Therefore,

\begin{eqnarray}
	S^{*}(t) \geq \mathrm{exp}^{(-\frac{\beta c}{N}-\mu)t}>0  
\end{eqnarray}
The positivity of the other components on $[0, \delta]$ can be shown in a similar manner. By the continuity of $S^{*}$ on $[0, \xi), (8)$ holds on $[0, \xi) .$ Therefore, $S^{*}(t)$ is bounded from below by
a positive number on $[0, \xi) .$ The same fact holds for $E^{*}$, $I_S^{*}$, $I_N^{*}$ and $T^{*}$ on $[0, \xi) .$ We indeed have,
$$\mathrm{lim} _{t \rightarrow \xi^{-}}(S^{*}(t), E^{*}(t), I_S^{*}(t), I_N^{*}(t), T^{*}(t))>0.$$ Therefore, $(S^{*}(t), E^{*}(t), I_S^{*}(t), I_N^{*}(t), T^{*}(t))$ can be continuously extended to $[0, \xi]$ and $\mathbb{R}_{+}^{5}$ is invariant with respect to the flow induced by system $(1)-(5)$ on $[0, \xi] .$ Consequently, $[0, \infty)$ is the maximal interval of existence of $\left(S^{*}(t), E^{*}(t), I_S^{*}(t), I_N^{*}(t), T^{*}(t)\right) .$ Therefore, $\mathbb{R}_{+}^{5}$ is positively invariant. \qed

\begin{theorem}
	For all the disease-free equilibrium points of the system of equations (1)-(5), we have $E=I_S=I_N=T=0$ and $S=\frac{\Lambda}{\mu}$.
\end{theorem}

\noindent{\textbf{Proof:}} Let $(\Tilde{S},\Tilde{E}, \Tilde{I_S}, \Tilde{I_N},\Tilde{T})$ be a disease-free equilibrium point. Then $\Tilde{E}=\Tilde{I_S}=\Tilde{I_N}=0$ and $\frac{d\Tilde{S}}{dt}=\frac{d\Tilde{E}}{dt}=\frac{d\Tilde{I_S}}{dt}=\frac{d\Tilde{I_N}}{dt}=\frac{d\Tilde{T}}{dt}=0$. Substituting these values to our system of equation (1)-(5), we have the following:
$$
\begin{array}{rcl}
0 & = & \mu \Tilde{T} \\
0 & = & \Lambda -\mu \Tilde{S}.
\end{array}
$$

Thus, $\Tilde{S}=\frac{\Lambda}{\mu}$. Since $\mu >0$, $\Tilde{T}=0$. \qed

$\;$\\
Since infected individuals are in $E$, $I_S$ and $I_N$, the rate of new infections in each compartment $(\mathcal{F})$ and the rate of other transitions between compartments
$\mathcal{V}$ can be rewritten as 
$$\mathcal{F}=\left(\begin{array}{c}
\beta cS\frac{I}{N}+\sigma \beta cT\frac{I}{N} \\
\alpha p\beta cE\frac{I}{N}\\
(1-\alpha)p\beta cE\frac{I}{N}
\end{array}\right),$$ \\ $$\mathcal{V}=\left(\begin{array}{c}
p\beta cE\frac{I}{N}+(\mu +k)E \\
-\alpha kE+ (\mu +r+d)I_S \\
-(1-\alpha)kE +(\mu +d)I_N
\end{array}\right).$$ Thus

$$F=\left(\begin{array}{ccc}
0 &\beta c & \beta c  \\
0 & 0 & 0\\
0 & 0 & 0
\end{array}\right),$$ \\
$$V=\left(\begin{array}{ccc}
\mu +k & 0 & 0 \\
-\alpha k & \mu +r+d & 0 \\
-(1-\alpha)k & 0 & \mu +d
\end{array}\right) and \quad V^{-1}=\left(\begin{array}{ccc}
\frac{1}{\mu +k} & 0 & 0\\
\frac{\alpha k}{(\mu +k)(\mu +r+d)} &\frac{1}{\mu +r+d} & 0\\
\frac{(1-\alpha)k}{(\mu +k)(\mu +d)} & 0 & \frac{1}{\mu +d}
\end{array}\right)$$

Therefore, the next generation matrix is

$$FV^{-1}=\left(\begin{array}{ccc}
\frac{\beta c\alpha k}{(\mu +k)(\mu +r+d)}+\frac{\beta c (1-\alpha)k}{(\mu +k)(\mu +d)} & \frac{\beta c}{\mu +r+d}& \frac{\beta c}{\mu +d} \\
0 & 0 & 0 \\
0 & 0 & 0
\end{array}\right).$$
Hence, the basic reproduction number is

$$R_0=\frac{\beta c\alpha k}{(\mu +k)(\mu +r+d)}+\frac{\beta c (1-\alpha)k}{(\mu +k)(\mu +d)}.$$\qed 

\begin{rem}
	In the model from Feng et.al.~\cite{13} the basic reproduction number corresponds to the case $\alpha=1$, since the authors did not consider the effect of stigmatization. 
\end{rem}

In the model of Feng et.al.\cite{13} there exist multiple endemic equilibria in the case of a strong reinfection. The same result is of course true for $\alpha=1$. We hence obtain as in \cite{13}

\begin{theorem}
	Let $\alpha=1$, i.e.~ $R_0=\frac{\beta c k}{(\mu +k)(\mu +r+d)}$. Then we have:
	\begin{itemize}
		\item[a)] If $R_0>1$, the system has exactly one endemic equilibrium point.
		\item[b)] If $R_0<1$ and $p>p_0$ for each $p$ there is an $R_p<1$, such that
		\begin{itemize}
			\item the system has two endemic equilibria if $R_0>R_p$
			\item the system has exactly one endemic equilibrium if $R_0=R_p$
			\item the system has no endemic equilibrium if $R_0<R_p$
		\end{itemize}
		\item[c)] If $R_0<1$ and $p<p_0$, the system has no endemic equilibrium point.	
		\item[d)] If $R_0<1$ and $p=p_0$, the system has one endemic equilibrium point.
	\end{itemize}	
\end{theorem}

For the proof see \cite{13}

For the analysis of endemic equilibria for $\alpha \neq 1$ assume $I^*>0$, where we use again $I=I_S+I_N$. 
We obtain then by setting the left hand side zero in (1)-(5) and with the abbreviation:
$$I_S^*=\alpha I^*$$
$$I_N^*=(1-\alpha) (1+\frac{r\alpha}{\mu+d})I^*$$
and hence 
$$(1-\alpha)\frac{r\alpha}{\mu+d} I^*=0$$
which can just be the case if $r=0$ or $\alpha=0$. In both cases the treatment plays no role. For $r=0$ one hence has the same behavior again as in \cite{13} for a vanishing treatment rate. It therefore is enough to consider $\alpha=0$. \\

We summarize this case in the following theorem.

\begin{theorem}
	Let $\alpha=0$, i.e.~ $R_0=\frac{\beta c k}{(\mu +k)(\mu+d)}$. Then we have:
	\begin{itemize}
		\item[a)] If $R_0>1$, the system has exactly one endemic equilibrium point.
		\item[b)] If $R_0<1$ and $p>p_0$ for each $p$ there is an $R_p<1$, such that
		\begin{itemize}
			\item the system has two endemic equilibria if $R_0>R_p$
			\item the system has exactly one endemic equilibrium if $R_0=R_p$
			\item the system has no endemic equilibrium if $R_0<R_p$
		\end{itemize}
		\item[c)] If $R_0<1$ and $p<p_0$, the system has no endemic equilibrium point.	
		\item[d)] If $R_0<1$ and $p=p_0$, the system has one endemic equilibrium point.
	\end{itemize}	
\end{theorem}

\noindent{\textbf{Proof:}}
Let $\alpha=0$, then we obtain 
$I_S^*=0$ and $I_N^*=I*$. Hence $T^*=0$.
As in \cite{13} we furthermore obtain 
$$S^*=\frac{\Lambda}{\mu + \beta c \frac{I^*}{N^*}}$$
$$E^*=\frac{\mu+d}{k+p\beta c \frac{I^*}{N^*}}I^*.$$

Now we use $N^*=\frac{\Lambda}{\mu}$ and the equation for the $E$ compartment and obtain with $x= \frac{I^*}{N^*}$:

$$
0 = \beta c S^* x - (p \beta c x +(\mu+k))E^* 
$$
Hence
$$0 = \beta c \frac{\Lambda}{\mu + \beta c x} x - (p \beta c x +(\mu+k))\frac{\mu+d}{k+p\beta c x}\frac{\Lambda}{\mu}x
$$
Getting rid of the $x$ in the denominator and dividing by $\frac{\Lambda}{\mu}x$ leads to: 
$$
0 = \mu \beta c (k+p \beta c x) - ((\mu+d)p \beta c x +(\mu+d)(\mu+k))(\mu + \beta c x)  
$$
Which gives
$$
0= x^2 + (1-\frac{1}{\mu +d} \frac{\mu}{\beta c} +\frac{\mu+k}{p \beta c})x - \frac{(\mu+k)\mu}{p \beta^2 c^2}(R_0-1).
$$
Due to feasibility the solution has to be real and positive. 

To obtain a real solution we need that
$$(1-\frac{1}{\mu +d} \frac{\mu}{\beta c} +\frac{\mu+k}{p \beta c})^2 + \frac{(\mu+k)\mu}{p \beta^2 c^2}(R_0-1)\geq 0$$

Hence:
$$
(p \beta c - \frac{\mu}{\mu+d}+ \mu +k)^2 + p (\mu +k) \mu (R_0-1)\geq0
$$
or 
$$
(p  - \frac{\mu}{(\mu+d)\beta c}+ \frac{\mu +k}{\beta c})^2 + p \frac{(\mu +k) \mu}{\beta^2 c^2} (R_0-1)\geq 0
$$
which is  always the case if  $R_0 >1$.

Let us hence take a look at the case $R_0\leq 1$, then we have an inequality of the form
$$
(p-a)^2 +pb \geq0
$$
which leads to
$$
p > \sqrt{\frac{b^2}{4}-ab} -a +\frac{b}{2} 
$$
Hence in our case 
\begin{multline*}
	p\geq \sqrt{\frac{(\mu +k)^2 \mu^2}{4\beta^4 c^4} (R_0-1)^2 - (\frac{\mu}{(\mu+d)\beta c}+ \frac{\mu +k}{\beta c})\frac{(\mu +k) \mu}{\beta^2 c^2} (R_0-1)}\\
	+ \frac{(\mu +k) \mu}{2\beta^2 c^2} (R_0-1) -\frac{\mu}{(\mu+d)\beta c}+ \frac{\mu +k}{\beta c} =p_0
\end{multline*}

Now let us check the different cases and how many solutions we obtain if we have real solution, i.e.~$p\geq p_0$: 

We consider the quadratic equation:
$x^2 +Px +Q =0$ with $P= 1-\frac{1}{\mu +d} \frac{\mu}{\beta c} +\frac{\mu+k}{p \beta c}$ and $Q= -\frac{(\mu+k)\mu}{p\beta^2 c^2}(R_0-1)$.

The solution is giving by 

$$
x_{1,2}= -\frac{P}{2}\pm \sqrt{\frac{P^2}{4}-Q}
$$

We have one (feasible) solution if either the term under the root is zero or there exist a positive and a negative solution, where the latter is to neglect.

We have $$\frac{P^2}{4}-Q=0 \Leftrightarrow (1-\frac{1}{\mu +d} \frac{\mu}{\beta c} +\frac{\mu+k}{p \beta c})^2 =- 4 \frac{(\mu+k)\mu}{p\beta^2 c^2}(R_0-1) $$

Hence we define
$$
R_p = 1- \frac{p\beta^2 c^2}{(\mu+k)\mu} (1-\frac{1}{\mu +d} \frac{\mu}{\beta c} +\frac{\mu+k}{p \beta c})^2.
$$
In that case if $R_0=R_p$ the term under the root vanishes. 

We have the following cases:

In the case $R_0>1$ and $p\geq p_0$ we have
$$x^2+Px+Q \text{ with } Q<0:$$ 
Hence we obtain just one positive solution and hence we have one endemic equilibrium.

In the case $R_0<1$ we have $Q>0$ which can lead to 3 different cases:

Now if $R_P>R_0$ the term under the root turns negative and there is no feasible solution. 

If  $R_P=R_0$ we have one solution hence one endemic equilibrium.

If $R_P<R_0$, $P<0$ and $p>p_0$ we have two solutions hence two endemic equilibria.

\qed

The stability analysis is a complete analogue of the proofs in \cite{13}.

\section{Optimal Control}
\subsection{Choice of controls}
We consider four control strategies where the first two, are controls to combat stigmatization. 

First, is the control minimizing the proportion of individuals going to the $I_N$ compartment from $E$, denoted by $u_1(t)$. Second is the control of encouraging infected people who are unwilling to get treated (in $I_N$) to change their views (to $I_S$), denoted by $u_2(t)$. These anti-stigmatization controls can be possibly done by advertising the positive effect of getting medically treated and downplaying the negative social connotation of being a TB patient. Another way of doing this is by directly giving money to infected people to support them and their families during the treatment. Although this method could be quite costly for the government, it may prove very effective in convincing the poor patients to seek and finish medical treatment.

The third control is increasing the treatment rate denoted by $u_3(t)$. This could be done by increasing the budget for tuberculosis treatment. However, because of the recent pandemic, budgets intended for TB treatment are cut to give more fund to combating the pandemic. This is the case specially in the developing countries. The fourth and last control minimizes the reinfection from the treated ($T$) compartment, denoted by $u_4(t)$. This control includes efforts of shielding the treated population against re-exposure from TB and their own efforts to increase their immune system.

For all time $t\geq 0$, the controls are in the interval $[0,1]$, with $0$ means a control is not implemented at all, and $1$ means full implementation of a control.

Our system with the controls is given by
\begin{eqnarray*}
	\frac{d S}{d t} &=\Lambda-\beta cS\frac{I}{N}-\mu S \\
	\frac{d E}{d t} &= \beta cS\frac{I}{N}-p\beta cE\frac{I}{N}-(\mu +k)E+(1-u_4(t))\sigma \beta c T\frac{I}{N} \\
	\frac{d I_S}{d t} &= (1+u_1(t))\alpha p\beta c E \frac{I}{N}+(1+u_1(t))\alpha kE-(\mu + (1+u_3(t))r+d)I_S + u_2(t)I_N\\
	\frac{d I_N}{d t} &= (1-(1+u_1(t))\alpha)p\beta cE\frac{I}{N} +(1-(1+u_1(t))\alpha)kE-(\mu +d)I_N - u_2(t)I_N\\
	\frac{d T}{d t} &= (1+u_3(t))rI_S-(1-u_4(t))\sigma \beta cT\frac{I}{N}-\mu T
\end{eqnarray*}

\subsection{Pontryagin's Maximum Principle}
We aim to minimize the number of exposed and infected individuals with the minimum implementation cost of the control measures. The objective function to be minimized is given by
\begin{eqnarray*}
	J(u_1,u_2,u_3,u_4)=\displaystyle\int_{t_0}^{t_f} \left(E(t) + I_N(t) + I_S(t) + \sum_{i=1}^4\frac{C_i}{2}u_i^2(t) \right) dt
\end{eqnarray*}

and the corresponding Hamiltonian $H$ is given by
\begin{eqnarray*}
	H=\displaystyle E(t) + I_N(t) + I_S(t) + \sum_{i=1}^4\frac{C_i}{2}u_i^2(t)+\sum_{i=1}^6\lambda_ig_i
\end{eqnarray*}
where $g_i$ is the right hand side of the differential equation of the $ith$ state variable. It is assumed that the controls are quadratic functions to incorporate nonlinear societal cost associated with the implementation of control measures. 

Applying Pontryagin's Maximum Principle, there exist adjoint variables $\lambda_1,...,\lambda_6$ which satisfy the following system of ordinary differential equations

\begin{eqnarray*}
	\frac{\partial\lambda_1}{\partial t}=&\lambda_1 \Bigg(\frac{\beta cI(N-S)}{N^2}+\mu\Bigg)-\lambda_2 \Bigg (\frac{\beta cI(N-S)+p\beta cEI-(1-u_4(t))\sigma\beta cTI}{N^2}\bigg )\\
	&+\lambda_3 \Bigg (\frac{(1+u_1(t))\alpha p\beta cEI}{N^2}\Bigg ) +\lambda_4 \Bigg (\frac{(1-(1+u_1(t))\alpha) p\beta cEI}{N^2}\Bigg )\\
	&-\lambda_5 \Bigg (\frac{(1-u_4(t))\sigma \beta cTI}{N^2}\Bigg )+\lambda_6\mu
\end{eqnarray*}

\begin{eqnarray*}
	\frac{\partial\lambda_2}{\partial t}=&-1-\lambda_1\Bigg(\frac{\beta cSI}{N^2}\Bigg ) + \lambda_2\Bigg (\frac{\beta cI(S+(1-u_4(t))\sigma T+p\beta cI(N-E)}{N^2}+u+k\Bigg )\\
	& -\lambda_3 \Bigg(\frac{(1+u_1(t))\alpha p\beta cI (N-E)}{N^2}+(1+u_1(t))\alpha k\Bigg )\\ &-\lambda_4 \Bigg[(1-(1+u_1(t))\alpha)\Bigg (k+\frac{p\beta cI(N-E)}{N^2}\Bigg )\Bigg]- \lambda_5 \Bigg (\frac{(1-u_4(t))\sigma \beta cTI}{N^2}\Bigg ) +\lambda_6 \mu
\end{eqnarray*}

\begin{eqnarray*}
	\frac{\partial\lambda_3}{\partial t}=&-1+\lambda_1\Bigg (\frac{\beta cS(N-I)}{N^2}\Bigg ) - \lambda_2 \Bigg (\frac{(\beta cS-p\beta cE+(1-u_4(t))\sigma\beta cT(N-I)}{N^2}\Bigg ) \\
	& -\lambda_3\Bigg (\frac{(1+u_1(t))\alpha p\beta cE(N-I)}{N^2}-(u+(1+u_3(t))r+d)\Bigg ) \\
	& -\lambda_4 \Bigg (\frac{(1-(1+u_1(t))\alpha)p\beta cE(N-I)}{N^2}\Bigg ) \\
	& -\lambda_5 \Bigg ((1+u_3(t))r-\frac{(1-u_4(t))\sigma \beta cT(N-I)}{N^2}\Bigg )+\lambda_6 (\mu +d)
\end{eqnarray*}

\begin{eqnarray*}
	\frac{\partial \lambda_4}{\partial t}= & -1+\lambda_1\Bigg (\frac{\beta cS(N-I)}{N^2}\Bigg )-\lambda_2 \Bigg (\frac{(\beta cS-p\beta cE +(1-u_4(t))\sigma \beta c T)(N-I)}{N^2}\Bigg )\\
	& -\lambda_3 \Bigg (\frac{(1+u_1(t))\alpha p \beta cE (N-I)}{N^2}+u_2(t)\Bigg ) +\lambda_4\Bigg (\frac{(1-(1+u_1(t))\alpha )p\beta cE (I-N)}{N^2}+\mu+d+u_2(t)\Bigg )\\ & +\lambda_5 \Bigg (\frac{(1-u_4(t))\sigma \beta cT (N-I)}{n^2}\Bigg ) +\lambda_6 (\mu+d)
\end{eqnarray*}

\begin{eqnarray*}
	\frac{\partial\lambda_5}{\partial t}= &-\lambda_1 \Bigg (\frac{\beta cSI}{N^2}\Bigg )-\lambda_2 \Bigg (\frac{-\beta cSI+p\beta cEI+(1-u_4(t))\sigma\beta cI (N-T)}{N^2} \Bigg )+\lambda_3 \Bigg (\frac{(1+u_1(t))\alpha p\beta c EI}{N^2}\Bigg )\\
	& +\lambda_4 \Bigg (\frac{(1-(1+u_1(t))\alpha )p\beta cEI}{N^2}\Bigg ) +\lambda_5 \Bigg (\frac{(1-u_4(t))\sigma\beta cI(N-T)}{N^2}+\mu\Bigg )+\lambda_6\mu
\end{eqnarray*}

\begin{eqnarray*}
	\frac{\partial\lambda_6}{\partial t}= &- \lambda_1 \Bigg (\frac{\beta cSI}{N^2}\Bigg )+ \lambda_2 \Bigg (\frac{\beta cSI-p\beta cEI+(1-u_4(t))\sigma\beta cTI)}{N^2}\Bigg )+\lambda_3 \Bigg(\frac{(1+u_1(t))\alpha p\beta cEI}{N^2} \Bigg ) \\
	& +\lambda_4 \Bigg (\frac{(1-(1+u_1(t))\alpha)p\beta cEI}{N^2}\Bigg )-\lambda_5 \Bigg(\frac{(1-u_4(t))\sigma \beta cTI}{N^2} \Bigg ) +\lambda_6\mu 
\end{eqnarray*}

The following theorem shows the form of our controls.

\begin{theorem}
	The optimal control variables are given by
	\begin{eqnarray*}
		u_1(t)= & \mathrm{min} \Bigg(1, \mathrm{max}\Bigg(0,\dfrac{(\lambda_4-\lambda_3)\Big (\alpha p\beta cE\frac{I}{N}+\alpha kE\Big)}{C_1}\Bigg)\\
		u_2(t)= & \mathrm{min} \Bigg(1, \mathrm{max}\Bigg(0,\frac{(\lambda_4-\lambda_3)I_N}{C_2}\Bigg) \Bigg)\\
		u_3(t)= & \mathrm{min} \Bigg(1, \mathrm{max}\Bigg(0,\frac{(\lambda_3-\lambda_5)rI_S}{C_3}\Bigg) \Bigg)\\
		u_4(t)= & \mathrm{min} \Bigg(1, \mathrm{max}\Bigg(0,\frac{(\lambda_2-\lambda_5) \sigma \beta cTI}{NC_4}\Bigg) \Bigg).
	\end{eqnarray*}
\end{theorem}

\noindent{\textbf{Proof:}}
Optimal controls $u_1,u_2,u_3$ and $u_4$ are derived by the following conditions:
\begin{eqnarray*}
	\frac{\partial H}{\partial u_1(t)}= & C_1u_1(t)+(\lambda_3-\lambda_4) \Big (\alpha p\beta cE\frac{I}{N}+\alpha kE\Big )=0 \\
	\frac{\partial H}{\partial u_2(t)}= & C_2u_2(t)+(\lambda_3-\lambda_4)I_N=0 \\
	\frac{\partial H}{\partial u_3(t)}= & C_3u_3(t)+(\lambda_5-\lambda_3)rI_S=0\\
	\frac{\partial H}{\partial u_4(t)}= & C_4u_4(t)+(\lambda_5-\lambda_2)\Big (\sigma\beta cT\frac{I}{N}\Big )=0
\end{eqnarray*} \qed

\section{Simulations}
\subsection{Parameter Values}
The values of the parameters used in the following simulations are given in Table~\ref{tab:parameters}. The values for the parameters $\mu$, $k$, $d$, and $r$ are from \cite{12}, while the values for $\Lambda$, $\sigma$, and $p$ are from \cite{13}. We also based on \cite{13} our choice for the initial conditions: $S(0)=18000$, $E(0)=5500$, $I_S(0)=700$, $I_N(0)=400$, and $T(0)=400$. The parameters $\beta$ and $c$ are separate parameters in \cite{13}, denoting the average numbers of susceptible infected by one infectious individual per contact per unit of time and the per-capita contact rate, respectively. But one may interpret the product $\beta c$ as just the transmission rate from $S$ to $E$. Its value in Table \ref{tab:parameters} is just an estimate to have an $R_0$ of around $3$ for the base model. In our simulations, we vary the level of stigmatization. We choose the values $\alpha = 0.3, 0.5, 0.7$ to represent high, medium, and low levels of stigmatization, respectively. For the weights of our controls $C_i, i=1,2,3,4$, we use the values $10$, $10^2$, and $10^3$ to denote low, medium, and high cost, respectively. Moreover, in the optimal control simulations we have the following lower and upper bounds for the controls: $0.01\leq u_1 \leq \frac{1-\alpha}{\alpha}$, $0.01\leq u_2 \leq 0.9$, $0.01\leq u_3 \leq \frac{1-r}{r}$, and $0.01\leq u_4 \leq 0.9$. We let the controls $u_1$ and $u_3$ increase the values of the parameters $\alpha$ and $r$ up to twice its given values but not greater than $1$.

\subsection{The Effect of Stigmatization}\label{effect}
\begin{figure}[htbp]
	\centering
	\includegraphics[scale=0.37]{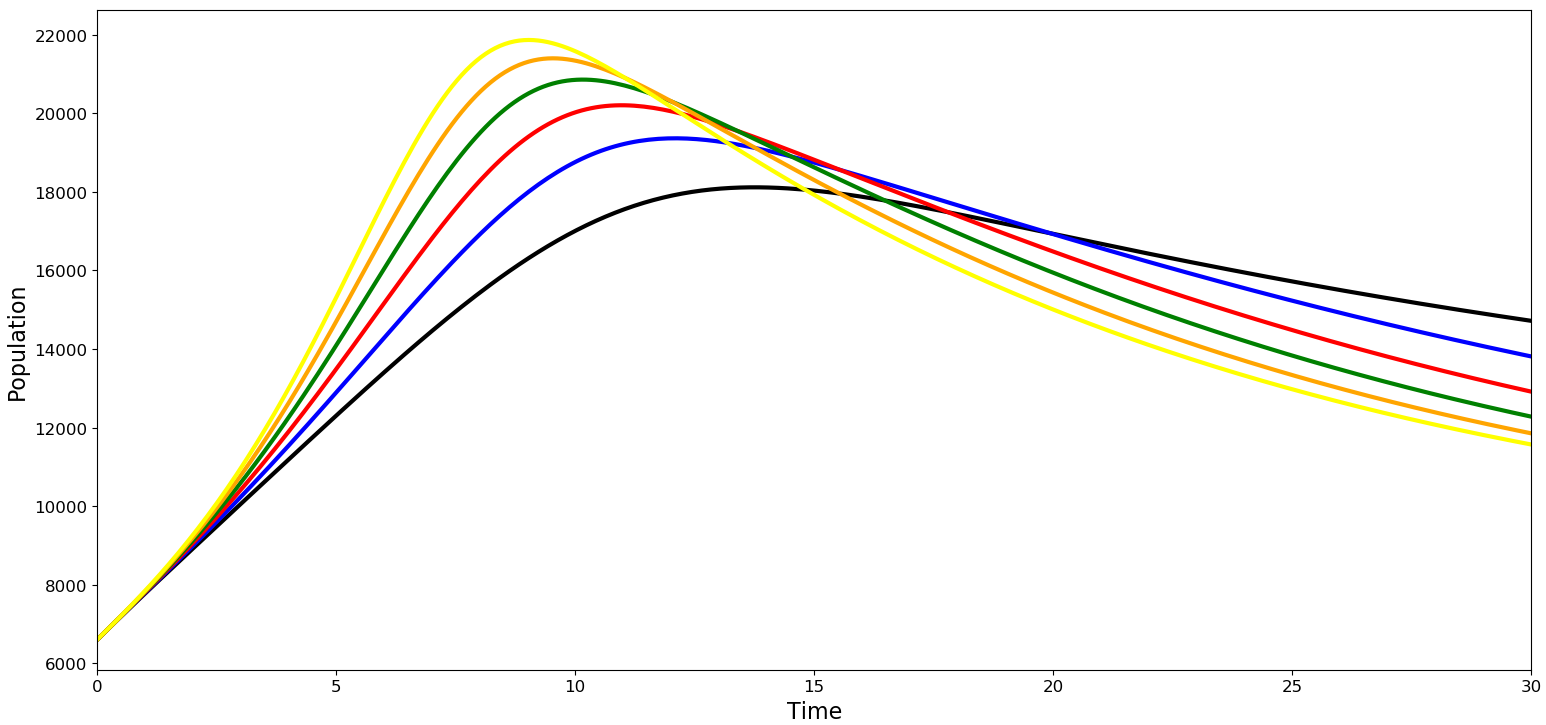} 
	\\[-0.5cm]\caption{The total infected $(E+I_S+I_N)$ over time. The curves with colors black, blue, red, green, orange, and yellow are for the simulations with values of $\alpha$ equal to $1$, $0.8$, $0.6$, $0.4$, $0.2$, and $0$, respectively.}
	
	\label{fig0}
\end{figure}

Using the model (\ref{eq1})-(\ref{eq5}), we simulate the effect of stigmatization by varying the values of $\alpha$. The results are given in Figure \ref{fig0} and Table \ref{tab0}.

\begin{table}[htbp]
	\centering
	\includegraphics[scale=0.45]{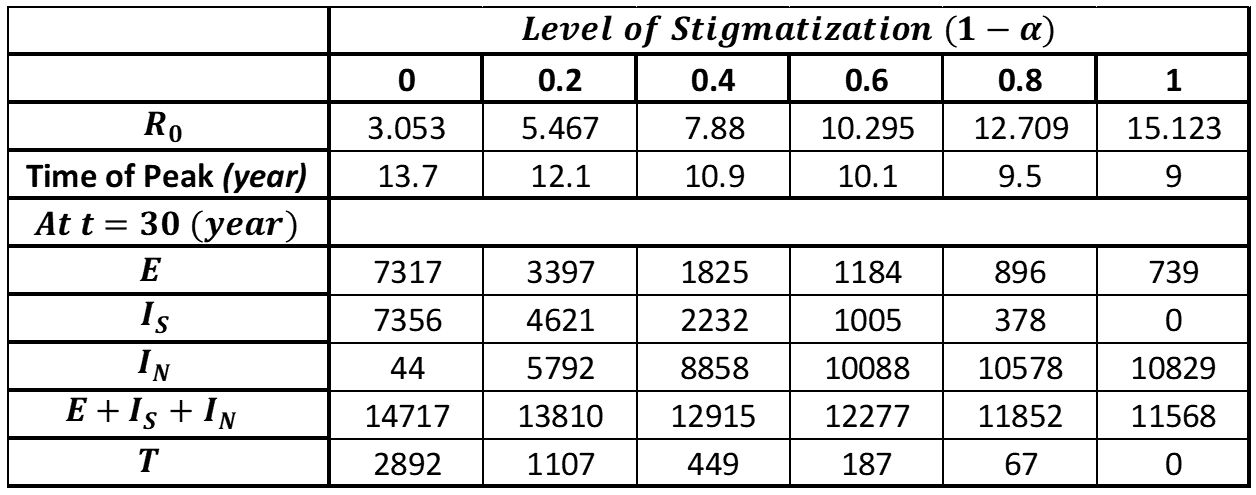} 
	\\[-0.5cm]\caption{The effect of the various levels of stigmatization in the dynamics of the TB model.}
	\label{tab0}
\end{table}
\subsection{Optimal Controls}\label{optimal}
\begin{figure}[h]
		\centering
		\includegraphics[width=\textwidth]{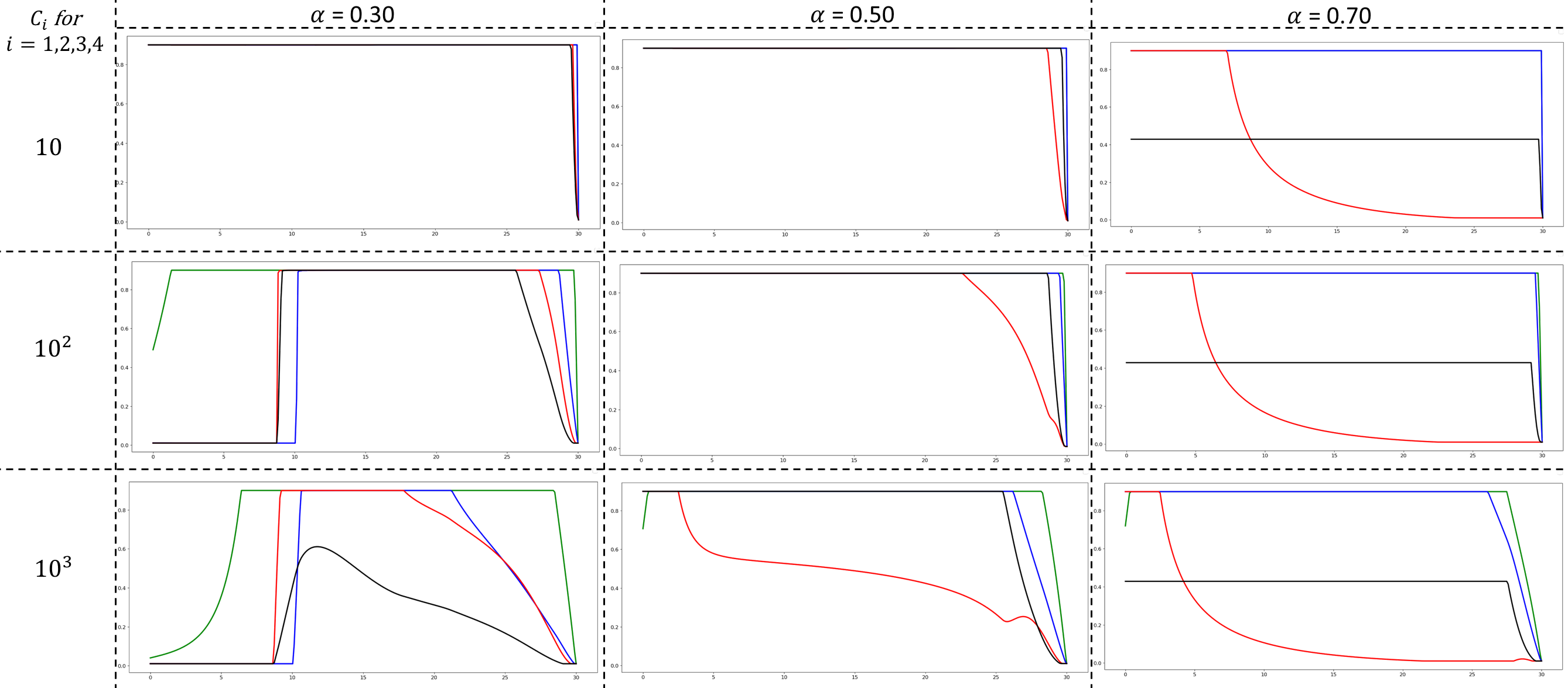} 
		\caption{Optimal controls. The $y$-axis represents the control rate from $0$ to $1$ and the $x$-axis represents the time in years. The simulations are up to 30 years. The curves black, red, blue, and green, are for the controls $u_1$, $u_2$, $u_3$, and $u_4$, respectively}
		\label{fig1}
	\end{figure}

In these simulations, we seek the optimal controls $(u_1, u_2, u_3, u_4)$ considering the cost of the controls $(C_i, i=1,2,3,4)$ and the level of stigmatization $(\alpha)$. The results are given in Figure \ref{fig1} and Table \ref{tab1}.

\begin{table}[htbp]
	\centering
	\includegraphics[scale=0.48]{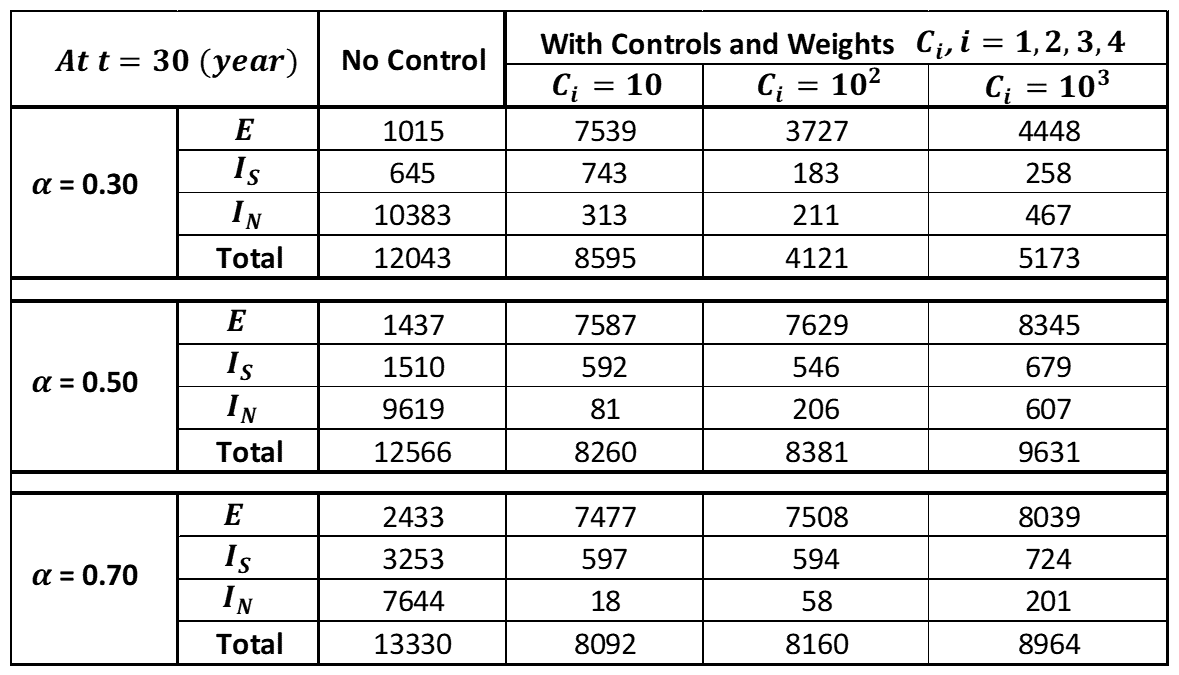} 
	\caption{The values of the infected compartments at $t=30$ for the various combination of control cost and level of stigmatization.}
	\label{tab1}
\end{table}

\subsection{Stigma Controls vs Treatment and Reinfection Controls}\label{stigmavsothers}
\begin{figure}[htbp]
	\centering
	\includegraphics[scale=0.37]{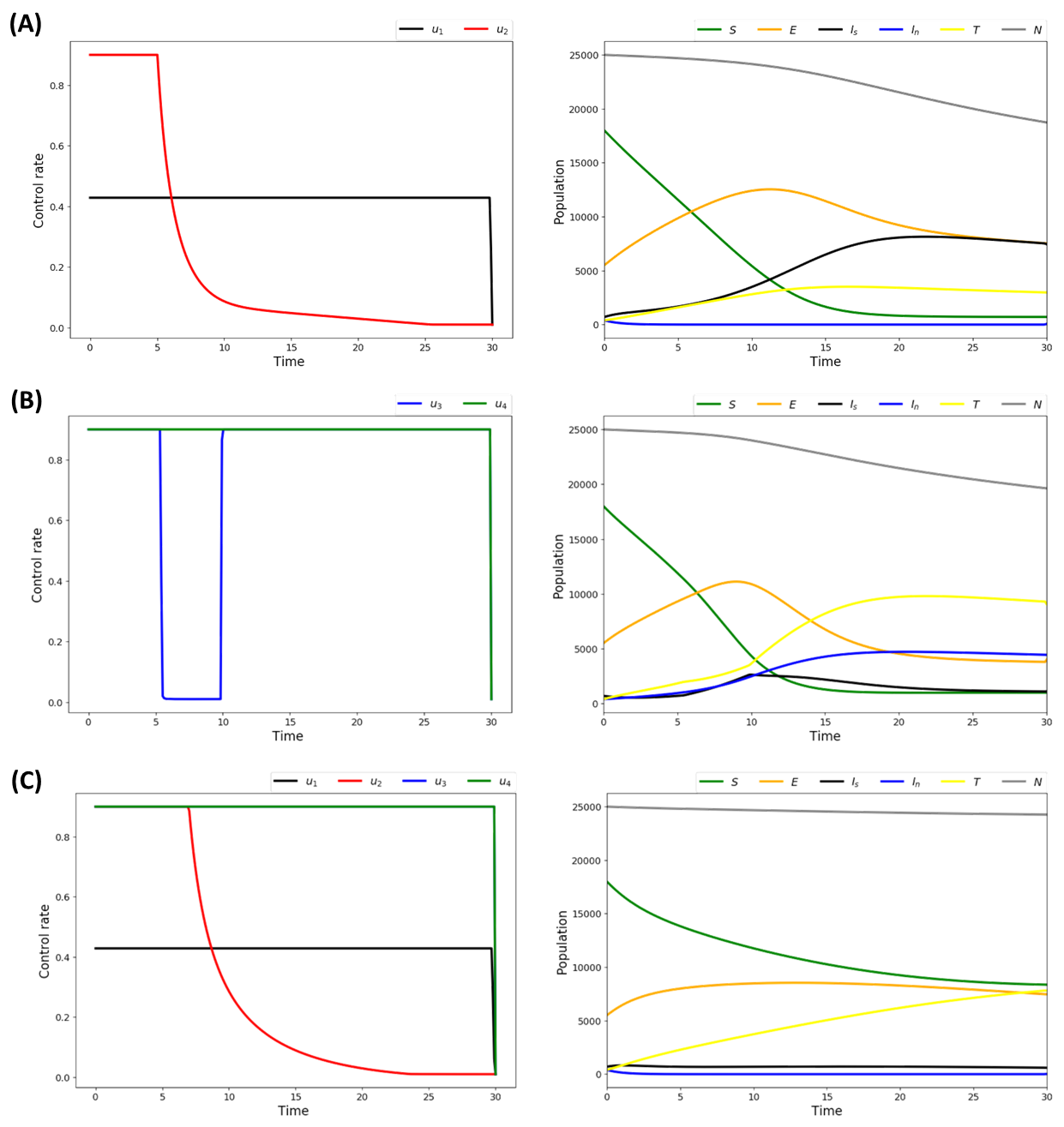} 
	\caption{For (A), (B), and (C), the graphs on the right are for the controls while the graphs on the left are for the compartments over time. In (A), we simulate having the stigma controls only, while in (B) having the other controls only. In (C), we simulate having all of the controls.}
	\label{fig2}
\end{figure}

We want to answer if the stigma controls ($u_1$ and $u_2$) are enough to curb the transmission of TB and if the controls are better compared to the other controls (treatment control $u_3$ and reinfection control $u_4$). In all the simulations, we use $C_i=10, i=1,2,3,4$ (low cost) and $\alpha = 0.7$ (low level of stigmatization). The results are given in Figure \ref{fig2} and Table \ref{tab2}.

\begin{table}[htbp]
	\centering
	\includegraphics[scale=0.45]{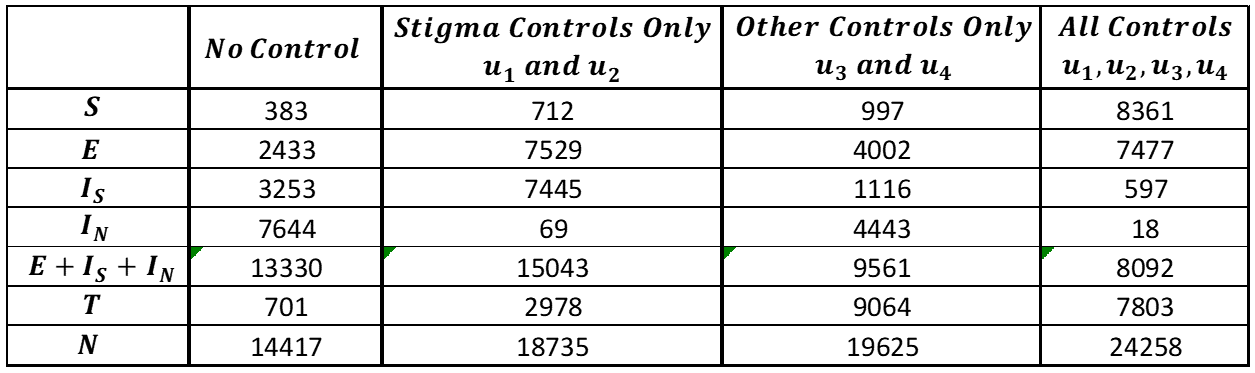} 
	\caption{Values of the compartments at $t=30$.}
	\label{tab2}
\end{table}

\section{Discussion}
From Section \ref{effect}, we conclude that as the level of stigmatization increases, the reproduction number increases considerably showing the negative effect of stigmatization. Although the total infected $(E+I_S+I_N)$ appears to be smaller for a higher level of stigmatization at the end of the simulation ($t=30$), one should notice that this is due to the fact that the for a higher $R_0$ the time to reach the peak is shorter. Moreover, one recognizes the effect of minimizing stigmatization in the number of treated individuals. We note that the more treated individuals, the more lives saved.

From the results in Section \ref{optimal}, we can see the various optimal controls for a particular level of cost and level of stigmatization. One could notice that in the optimal controls given in Figure \ref{fig1}, the control $u_4$ minimizing the reinfection from the treated compartment takes priority seconded by the control $u_3$ increasing the treatment rate. In Section \ref{stigmavsothers}, where we compare the relative impact of the treatment and reinfection controls $(u_3, u_4)$ with the stigma controls $(u_1, u_2)$. In Table \ref{tab2}, we can see clearly that the treatment and reinfection controls provide better results than the stigma controls. But the result is improved if all the controls are used together.

\section{Conclusion}
We developed and studied a tuberculosis model with exogenous reinfection and stigmatization. As in \cite{13} we elaborated the existence of multiple endemic equilibrium points in dependence of the reinfection parameter. The stigmatization leads to a more involved dynamics in the case where $\alpha \notin \{0,1\}$, in these cases one has no endemic equilibrium point. It is to expect that in that case there can be oscillations within the infected compartments. As in \cite{13} the basic reproductive number was calculated. In the case $\alpha=1$ we obtain the same result as \cite{Feng}. 
Our simulations showed the considerable negative effect of stigmatization. However, the optimal control results showed also that strategies relying only on minimizing stigmatization is far from enough in curbing the disease. In fact, the simulations show that the treatment and reinfection controls have better results if they are to be used separately. A better strategy is to use the four controls together and some optimal ways of doing it is given in Figure \ref{fig1}. One may notice in Table \ref{tab2}, last row, that having the stigma controls in addition to the other controls could mean more saved lives as can be seen in the increasing total population at the end time of the simulations $(t=30)$.

\newpage
\bibliographystyle{plain}



\end{document}